\documentclass[12pt]{article}
\usepackage{amssymb}
\usepackage{amsthm}
\newtheorem{Definition}{Definition}
\newtheorem{Lemma}[Definition]{Lemma}

\newtheorem{Theorem}[Definition]{Theorem}
\newtheorem{Remark}[Definition]{Remark}

\newtheorem{Corollary}[Definition]{Corollary}
\title{When does a semiring become a residuated lattice?\thanks{%
Preprint of an article published in Asian-European Journal of Mathematics, Vol.~9, No.~4, 2016, pages 1650088-1--1650088-10,
\texttt{DOI: 10.1142/S1793557116500881} \textcopyright\ World Scientific Publishing Company,
\mbox{\texttt{https://www.worldscientific.com/worldscinet/aejm}}%
}}
\author{Ivan~Chajda and Helmut~L\"anger}
\date{}
\begin{document}
\footnotetext[1]{Support of the research by the bilateral project "New perspectives on residuated posets", supported by the Austrian Science Fund (FWF), project I~1923-N25, and the Czech Science Foundation (GA\v CR), project 15-34697L, as well as by AKTION Austria -- Czech Republic, project 75p11, is gratefully acknowledged.}
\maketitle
\begin{abstract}
It is an easy observation that every residuated lattice is in fact a semiring because multiplication distributes over join and the other axioms of a semiring are satisfied trivially. This semiring is commutative, idempotent and simple. The natural question arises if the converse assertion is also true. We show that the conversion is possible provided the given semiring is, moreover, completely distributive. We characterize semirings associated to complete residuated lattices satisfying the double negation law where the assumption of complete distributivity can be omitted. A similar result is obtained for idempotent residuated lattices.
\end{abstract}



It is well-known that semirings are useful structures with applications both in mathematics and in computer science, see e.g.\ \cite{KS} for details. In particular, unitary rings as well as bounded distributive lattices are examples of semirings. We recall the definition of a semiring from the monograph by J.~S.~Golan (\cite G).

\begin{Definition}\label{def1}
A {\em semiring} is an algebra ${\mathcal S}=(S,+,\cdot,0,1)$ of type $(2,2,0,0)$ satisfying the following conditions for all $x,y,z\in S$:
\begin{enumerate}
\item[{\rm(i)}] $(S,+,0)$ is a commutative monoid,
\item[{\rm(ii)}] $(S,\cdot,1)$ is a monoid,
\item[{\rm(iii)}] $(x+y)\cdot z=x\cdot z+y\cdot z$ and $x\cdot(y+z)=x\cdot y+x\cdot z$,
\item[{\rm(iv)}] $x\cdot0=0\cdot x=0$.
\end{enumerate}
${\mathcal S}$ is called
\begin{itemize}
\item {\em idempotent} if it satisfies the identity $x+x\approx x$ {\rm(}i.e., $(S,+)$ is a semilattice{\rm)},
\item {\em commutative} if it satisfies the identity $x\cdot y\approx y\cdot x$,
\item {\em simple} {\rm(}see e.g.\ {\rm\cite G)} if it satisfies the identity $x+1\approx1$.
\end{itemize}
\end{Definition}

Recall that a {\em ring} ${\mathcal R}=(R,+,\cdot)$ is called {\em Boolean} if it satisfies the identity $x\cdot x=x$. Note that any Boolean ring is commutative and of characteristic $2$. If ${\mathcal R}$ is a Boolean ring with $1$ then it can be organized into a Boolean algebra. Also conversely, every Boolean algebra can be organized into a Boolean ring with $1$. Since Boolean algebras form an algebraic axiomatization of the classical propositional logic, we can ask if a similar role is played by semirings for certain non-classical logics. The answer is partly positive. The aim of this paper is to study this question.

For the concept of a residuated lattice we use the definition given in \cite B.

\begin{Definition}\label{def2}
A {\em residuated lattice} is an algebra ${\mathcal L}=(L,\vee,\wedge,\odot,\rightarrow,0,1)$ of type $(2,2,2,$ $2,0,0)$ satisfying the following conditions for all $x,y,z\in L$:
\begin{enumerate}
\item[{\rm(i)}] $(L,\vee,\wedge,0,1)$ is a bounded lattice,
\item[{\rm(ii)}] $(L,\odot,1)$ is a commutative monoid,
\item[{\rm(iii)}] $x\odot y\leq z$ if and only if $x\leq y\rightarrow z$.
\end{enumerate}
Here $\leq$ denotes the lattice order. Condition {\rm(iii)} is called the {\em adjointness property}. The operation $\rightarrow$ is called {\em residuum} and for $a,b\in L$, $a\rightarrow b$ is called the residuum of $b$ by $a$.
\end{Definition}

We recall some well-known facts concerning residuated lattices:

\begin{Lemma}\label{lem3}
Let ${\mathcal L}=(L,\vee,\wedge,\odot,\rightarrow,0,1)$ be a residuated lattice and $a,b\in L$. Then {\rm(i)} -- {\rm(vi)} hold:
\begin{enumerate}
\item[{\rm(i)}] $a\rightarrow b$ is the greatest element of $\{x\in L\,|\,x\odot a\leq b\}$,
\item[{\rm(ii)}] $\neg0=1$,
\item[{\rm(iii)}] $\neg1=0$,
\item[{\rm(iv)}] $a\odot\neg a=0$,
\item[{\rm(v)}] $a\leq\neg(\neg a)$,
\item[{\rm(vi)}] $a\leq b$ implies $\neg b\leq\neg a$.
\end{enumerate}
Here and in the following $\neg a$ is an abbreviation for $a\rightarrow0$ and is called the {\em negation} of $a$. Note that $\neg a$ is the greatest element of $\{x\in L\,|\,x\odot a=0\}$.
\end{Lemma}

\begin{proof}
Follows easily from the definition. See e.g.\ Theorems~2.17 and 2.30 of \cite B.
\end{proof}

Recall that residuated lattices form an algebraic axiomatization of fuzzy logics, see \cite B for details and for the complete explanation. Hence, when studying connections between semirings and residuated lattices, we in fact ask for an algebraic axiomatization of fuzzy logics by means of certain semirings in the way known for classical logic and Boolean rings.

\begin{Theorem}\label{th1}
Let ${\mathcal L}=(L,\vee,\wedge,\odot,\rightarrow,0,1)$ be a residuated lattice. Then
\[
{\bf S}({\mathcal L}):=(L,\vee,\odot,0,1)
\]
is a commutative idempotent simple semiring.
\end{Theorem}

\begin{proof}
Let $a,b,c\in L$. Since $(L,\vee,\wedge,0,1)$ is a bounded lattice, its reduct $(L,\vee,0)$ is a commutative monoid. Of course, $\vee$ is idempotent. Further, $(L,\odot,1)$ is a commutative monoid by definition. Since $0$ is the least element of ${\mathcal L}$, we have $0\leq a\rightarrow0$ and according to the adjointness property $0\odot a\leq0$ and thus $0\odot a=0$. The law
\[
a\odot(b\vee c)=(a\odot b)\vee(a\odot c)
\]
holds according to Theorem~2.25 of \cite B. Finally, $a\vee1=1$. Altogether, ${\bf S}({\mathcal L})=(L,\vee,\odot,0,1)$ is a commutative idempotent simple semiring.
\end{proof}

\begin{Remark}
Theorem~2.25 of {\rm\cite B} reads as follows: \\
{\rm(}distributivity of $\odot,\rightarrow$ over $\wedge,\vee${\rm)} The following are true for each index set $I$. Moreover, in the first three formulas, if the left side makes sense, then so does the right side.
\begin{eqnarray*}
x\odot\bigvee_{i\in I}y_i & = & \bigvee_{i\in I}(x\odot y_i), \\
x\rightarrow\bigwedge_{i\in I}y_i & = & \bigwedge_{i\in I}(x\rightarrow y_i), \\
\bigvee_{i\in I}x_i\rightarrow y & = & \bigwedge_{i\in I}(x_i\rightarrow y), \\
x\odot\bigwedge_{i\in I}y_i & \leq & \bigwedge_{i\in I}(x\odot y_i), \\
\bigvee_{i\in I}(x\rightarrow y_i) & \leq & x\rightarrow\bigvee_{i\in I}y_i, \\
\bigvee_{i\in I}(x_i\rightarrow y) & \leq & \bigwedge_{i\in I}x_i\rightarrow y.
\end{eqnarray*}
\end{Remark}

We want to investigate if the converse statement also holds. At first we need the following

\begin{Lemma}\label{lem1}
Let ${\mathcal S}=(S,+,\cdot,0,1)$ be an idempotent semiring and $\leq$ the induced order of the join-semilattice $(S,+)$. Then $\cdot$ is isotone, i.e., $a\leq b$ implies both $a\cdot c\leq b\cdot c$ and $c\cdot a\leq c\cdot b$, for all $a,b,c\in S$.
\end{Lemma}

\begin{proof}
Assume $a\leq b$. Using distributivity of $\cdot$ with respect to $+$ we have
\[
a\cdot c\leq a\cdot c+b\cdot c=(a+b)\cdot c=b\cdot c.
\]
The second inequality can be shown analogously.
\end{proof}

Let $(S,+)$ be a join-semilattice and $x_i\in S$ for all $i\in I$. Then the sum $\sum\limits_{i\in I}x_i$ denotes the supremum of the $x_i$ (if it exists). Of course, $\sum\limits_{i\in I}x_i$ exists in case $I$ is finite.

In order to obtain a residuated lattice from a given semiring we need the following concept.

\begin{Definition}\label{def3}
A commutative idempotent semiring ${\mathcal S}=(S,+,\cdot,0,1)$ is called {\em completely distributive} if the following conditions hold:
\begin{enumerate}
\item[{\rm(i)}] In the join-semilattice $(S,+)$, the sum $\sum\limits_{i\in I}y_i$ exists for all sets $I$ and $y_i\in S$ {\rm(}$i\in I${\rm)},
\item[{\rm(ii)}] $x\cdot(\sum\limits_{i\in I}y_i)=\sum\limits_{i\in I}(x\cdot y_i)$ for all $x,y_i\in S$ {\rm(}$i\in I${\rm)}.
\end{enumerate}
\end{Definition}

Note that both conditions hold in case $S$ is finite.

Now we can prove the following result.

\begin{Theorem}\label{th3}
Let ${\mathcal S}=(S,+,\cdot,0,1)$ be a commutative idempotent simple completely distributive semiring. Then there exist binary operations $\wedge$ and $\rightarrow$ on $S$ such that
\[
{\bf L}({\mathcal S}):=(S,+,\wedge,\cdot,\rightarrow,0,1)
\]
is a residuated lattice.
\end{Theorem}

\begin{proof}
Clearly, $(S,+,0)$ is a join-semilattice with least element $0$. Let $a,b,c\in S$. Because of the simplicity of ${\mathcal S}$ we have $a\leq1$ where $\leq$ denotes the induced semilattice order. In other words, $(S,+,0,1)$ is a bounded join-semilattice. According to (i) of Definition~\ref{def3}, in the poset $(S,\leq)$ there exist arbitrary suprema and hence also arbitrary infima. Therefore we obtain a bounded complete lattice $(S,+,\wedge,0,1)$. Put
\[
a\rightarrow b:=\sum_{a\cdot x\leq b}x.
\]
From this follows that $a\cdot b\leq c$ implies $a\leq b\rightarrow c$. Conversely, according to complete distributivity of ${\mathcal S}$ and isotonicity of $\cdot$, the assertion $a\leq b\rightarrow c$ implies
\[
a\cdot b\leq(b\rightarrow c)\cdot b=(\sum_{b\cdot x\leq c}x)\cdot b=\sum_{b\cdot x\leq c}(x\cdot b)\leq c.
\]
This completes the proof of the theorem.
\end{proof}

\begin{Corollary}\label{cor1}
Let ${\mathcal S}=(S,+,\cdot,0,1)$ be a finite commutative idempotent simple semiring. Then there exist binary operations $\wedge$ and $\rightarrow$ on $S$ such that
\[
{\bf L}({\mathcal S}):=(S,+,\wedge,\cdot,\rightarrow,0,1)
\]
is a residuated lattice.
\end{Corollary}

\begin{Remark}\label{rem1}
By the previous theorem we see that complete distributivity of commutative idempotent simple semirings implies the adjointness property and, conversely, the adjointness property implies distributivity of the derived semiring {\rm(}as shown e.g.\ in Theorem~2.25 in {\rm\cite B)}.
\end{Remark}

When studying residuated lattices, we usually ask for some additional conditions that are satisfied by algebras derived from residuated lattices. One of them is the so-called double negation law.

\begin{Definition}\label{def5}
A residuated lattice ${\mathcal L}=(L,\vee,\wedge,\odot,\rightarrow,0,1)$ is said to satisfy the {\em double negation law} if it satisfies the identity $\neg(\neg x)=x$.
\end{Definition}

For example, an MV-algebra is a residuated lattice satisfying the double negation law as well as the so-called divisibility, i.e.\ the identity $x\wedge y\approx x\odot(x\rightarrow y)$ (cf.\ Theorem~2.42 of \cite B).

For the reader's convenience, we recall the definition of an MV-algebra. An {\em MV-algebra} is an algebra $(L,\oplus,\neg,0)$ of type $(2,1,0)$ satisfying the following identities:
\begin{eqnarray*}
(x\oplus y)\oplus z & \approx & x\oplus(y\oplus z), \\
x\oplus y & \approx & y\oplus x, \\
x\oplus0 & \approx & x, \\
x\oplus\neg0 & \approx & \neg0, \\
\neg(\neg x) & \approx & x, \\
\neg(\neg x\oplus y)\oplus y & \approx & \neg(\neg y\oplus x)\oplus x.
\end{eqnarray*}
If one defines
\begin{eqnarray*}
x\rightarrow y & := & \neg x\oplus y, \\
x\vee y & := & (x\rightarrow y)\rightarrow y, \\
x\wedge y & := & \neg(\neg x\vee\neg y), \\
x\odot y & := & \neg(\neg x\oplus\neg y), \\
1 & := & \neg0
\end{eqnarray*}
then $(L,\vee,\wedge,\odot,\rightarrow,0,1)$ is a residuated lattice.

As pointed out in Theorem~\ref{th3}, if ${\mathcal S}=(S,+,\cdot,0,1)$ is a commutative idempotent simple semiring then $(S,+,0,1)$ is a bounded join-semilattice. Let again denote the induced order by $\leq$ and assume that every subset of $S$ has a supremum with respect to $\leq$. Then $(S,\leq)$ is a complete lattice and we can define a unary operation $n$ on $S$ as follows:
\begin{equation}\label{equ1}
n(x):=\sum_{x\cdot y=0}y\mbox{ for all }x\in S.
\end{equation}
The following result is an easy consequence of the definition of $n$.

\begin{Lemma}\label{lem2}
If ${\mathcal S}=(S,+,\cdot,0,1)$ is a commutative idempotent simple completely distributive semiring and $a\in S$ then $a\cdot n(a)=0$ and $a\leq n(n(a))$.
\end{Lemma}

\begin{proof}
We have
\[
a\cdot n(a)=a\cdot\sum_{a\cdot y=0}y=\sum_{a\cdot y=0}(a\cdot y)=0
\]
and hence
\[
n(n(a))=\sum_{n(a)\cdot y=0}y\geq a.
\]
\end{proof}

In what follows, we provide a correspondence between complete residuated lattices satisfying the double negation law and a certain class of semirings. Contrary to Theorem~\ref{th3} we do not need complete distributivity in this case.

\begin{Definition}\label{def4}
A {\em DNL-semiring} {\rm(double negation law semiring)} is a commutative idempotent simple semiring ${\mathcal S}=(S,+,\cdot,0,1)$ such that $(S,+)$ is a complete join-semilattice and the following identities are satisfied for the function $n$ defined by {\rm(\ref{equ1})}:
\begin{enumerate}
\item[{\rm(i)}] $n(n(x))\approx x$,
\item[{\rm(ii)}] $x\cdot n(x\cdot n(y))\leq y$ {\rm(}i.e., $x\cdot n(x\cdot n(y))+y\approx y${\rm)},
\item[{\rm(iii)}] $x\cdot n(x+y)\approx0$.
\end{enumerate}
\end{Definition}

We are now able to prove the following result:

\begin{Theorem}\label{th4}
Let ${\mathcal S}=(S,+,\cdot,0,1)$ be a {\rm DNL}-semiring and put
\begin{eqnarray*}
x\wedge y & := & n(n(x)+n(y))\mbox{ and} \\
x\rightarrow y & := & n(x\cdot n(y))
\end{eqnarray*}
for all $x,y\in S$. Then
\[
{\bf L}({\mathcal S}):=(S,+,\wedge,\cdot,\rightarrow,0,1)
\]
is a complete residuated lattice satisfying the double negation law.
\end{Theorem}

\begin{proof}
Let $a,b,c\in S$. Since ${\mathcal S}=(S,+,\cdot,0,1)$ is a DNL-semiring, $(S,\cdot,1)$ is a commutative monoid and $(S,\leq)$ a complete lattice with smallest element $0$ and greatest element $1$. Here $\leq$ denotes the induced order of ${\mathcal S}$. 
First of all we have to show that $n(n(x)+n(y))$ is the infimum of $x$ and $y$ with respect to $\leq$. According to (i) of Definition~\ref{def4}, $n$ is an involution. If $a\leq b$ then
\[
n(b)=\sum_{b\cdot y=0}y\leq\sum_{a\cdot y=0}y=n(a)
\]
since according to Lemma~\ref{lem1} we have for all $y\in S$ that $b\cdot y=0$ implies $a\cdot y=0$. Since $n$ is an involution we have $a\leq b$ if and only if $n(b)\leq n(a)$, i.e., $n$ is an order isomorphism from $(S,\leq)$ to $(S,\geq)$. Let $x\wedge y$ denote the infimum of $x$ and $y$ with respect to $\leq$. Then $n(x\wedge y)=n(x)+n(y)$ and hence $x\wedge y=n(n(x)+n(y))$. Next we want to prove the double negation law. Put
\[
\neg a:=a\rightarrow0.
\]
Then
\[
\neg a=n(a\cdot n(0))=n(a\cdot1)=n(a)
\]
and according to (i) of Definition~\ref{def4} we have $\neg(\neg a)=a$ proving the double negation law. Finally, we need to prove the adjointness property. If $a\leq b$ then according to (iii) of Definition~\ref{def4} we have
\[
a\rightarrow b=n(a\cdot n(b))=n(a\cdot n(a+b))=n(0)=1.
\]
The first equality holds according to the definition of $\rightarrow$. \\
Conversely, if $a\rightarrow b=1$ then $n(a\cdot n(b))=1$ and according to (ii) of Definition~\ref{def4}
\[
a=a\cdot1=a\cdot n(a\cdot n(b))\leq b.
\]
Hence
\[
a\leq b\mbox{ if and only if }a\rightarrow b=1.
\]
Using this equivalence and (i) of Definition~\ref{def4} we see that the following are equivalent:
\begin{eqnarray*}
a\cdot b & \leq & c, \\
a\cdot b\rightarrow c & = & 1, \\
n(a\cdot b\cdot n(c)) & = & 1, \\
n(a\cdot n(n(b\cdot n(c)))) & = & 1, \\
a\rightarrow n(b\cdot n(c)) & = & 1, \\
a & \leq & n(b\cdot n(c)), \\
a & \leq & b\rightarrow c.
\end{eqnarray*}
This shows the adjointness property. Altogether, ${\bf L}({\mathcal S})=(S,+,\wedge,\cdot,\rightarrow,0,1)$ is a complete residuated lattice satisfying the double negation law.
\end{proof}

\begin{Theorem}\label{th6}
Let ${\mathcal L}=(L,\vee,\wedge,\odot,\rightarrow,0,1)$ be a complete residuated lattice satisfying the double negation law. Then
\[
{\bf S}({\mathcal L}):=(L,\vee,\odot,0,1)
\]
is a {\rm DNL}-semiring.
\end{Theorem}

\begin{proof}
Let $a,b,c\in L$. According to Theorem~\ref{th1}, ${\bf S}({\mathcal L})$ is a commutative idempotent simple semiring. According to Theorem~2.40 in \cite B, $a\rightarrow b=\neg(a\odot\neg b)$ since ${\mathcal L}$ satisfies the double negation law. Hence, according to Theorem~2.17 in \cite B we have
\[
\neg a=a\rightarrow0=\bigvee_{a\cdot x=0}x=n(a).
\]
Thus $n(n(a))=\neg(\neg a)=a$ proving (i) of Definition~\ref{def4}. According to Theorem~2.18 in \cite B, every residuated lattice satisfies the following identities:
\begin{enumerate}
\item[(a)] $x\odot(x\rightarrow y)\leq y$ (i.e., $(x\odot(x\rightarrow y))\vee y\approx y$),
\item[(b)] $x\rightarrow(x\vee y)\approx1$.
\end{enumerate}
Identity (a) can be rewritten as
\[
x\odot\neg(x\odot\neg y)\leq y
\]
which is (ii) of Definition~\ref{def4}. Moreover, identity (b) can be rewritten as
\[
\neg(x\odot\neg(x\vee y))\approx1.
\]
Thus, according to the double negation law
\[
x\odot\neg(x\vee y)\approx0
\]
which is (iii) of Definition~\ref{def4}. Hence ${\bf S}({\mathcal L})$ is a DNL-semiring.
\end{proof}

We can prove that the correspondence described in the last two theorems is one-to-one.

\begin{Theorem}\label{th7}
$\mbox{}$
\begin{enumerate}
\item[{\rm(i)}] If ${\mathcal S}=(S,+,\cdot,0,1)$ is a {\rm DNL}-semiring then ${\bf S}({\bf L}({\mathcal S}))={\mathcal S}$.
\item[{\rm(ii)}] If ${\mathcal L}=(L,\vee,\wedge,\odot,\rightarrow,0,1)$ is a complete residuated lattice satisfying the double negation law then ${\bf L}({\bf S}({\mathcal L}))={\mathcal L}$.
\end{enumerate}
\end{Theorem}

\begin{proof}
$\mbox{}$
\begin{enumerate}
\item[(i)] If ${\bf L}({\mathcal S})=(S,+,\wedge,\cdot,\rightarrow,0,1)$ then ${\bf S}({\bf L}({\mathcal S}))=(S,+,\cdot,0,1)={\mathcal S}$.
\item[(ii)] We have ${\bf S}({\mathcal L})=(L,\vee,\odot,0,1)$. Put ${\bf L}({\bf S}({\mathcal L}))=(L,\vee,\cap,\odot,\Rightarrow,0,1)$ and let $a,b\in L$. Since the orders in ${\mathcal L}$ and ${\bf L}({\bf S}({\mathcal L}))$ coincide we have $a\cap b=a\wedge b$. Furthermore,
\begin{eqnarray*} \\
& & n(a)=\bigvee_{a\odot x=0}x=\bigvee_{a\odot x\leq0}x=a\rightarrow0=\neg a\mbox{ and} \\
& & a\Rightarrow b=n(a\odot n(b))=\neg(a\odot\neg b)=a\rightarrow b.
\end{eqnarray*}
Hence ${\bf L}({\bf S}({\mathcal L}))={\mathcal L}$.
\end{enumerate}
\end{proof}

Another useful property of residuated lattices is idempotency. A residuated lattice ${\mathcal L}=(L,\vee,\wedge,\odot,\rightarrow,0,1)$ is called {\em idempotent} if it satisfies the identity $x\odot x=x$.

\begin{Theorem}\label{th5}
A residuated lattice ${\mathcal L}=(L,\vee,\wedge,\odot,\rightarrow,0,1)$ is idempotent if and only if ${\bf S}({\mathcal L})$ is a bounded distributive lattice.
\end{Theorem}

\begin{proof}
Let $a,b\in L$. According to Theorem~2.17 in \cite B we have $a\odot b\leq a\wedge b$. Since the operation $\odot$ is monotone (Theorem~2.23 in \cite B) and $a\wedge b\leq a$ and $a\wedge b\leq b$, using idempotency we infer
\[
a\wedge b=(a\wedge b)\odot(a\wedge b)\leq a\odot b\leq a\wedge b
\]
proving that $\wedge$ and $\odot$ coincide. According to distributivity ((iii) of Definition~\ref{def1} and Theorem~\ref{th1}) we obtain that ${\bf S}({\mathcal L})$ is a bounded distributive lattice. The converse assertion is trivial.
\end{proof}

\begin{Corollary}\label{cor2}
If ${\mathcal L}$ additionally satisfies the double negation law then ${\mathcal L}$ is idempotent if and only if ${\bf S}({\mathcal L})$ is a Boolean algebra.
\end{Corollary}

\begin{proof}
This follows from the previous proof, from (iv) and (vi) of Lemma~\ref{lem3} and from the double negation law.
\end{proof}

\begin{Definition}\label{def6}
A residuated lattice $(L,\vee,\wedge,\odot,\rightarrow,0,1)$ is called {\em prelinear} if it satisfies the identity
\[
(x\rightarrow y)\vee(y\rightarrow x)\approx1.
\]
\end{Definition}

\begin{Corollary}\label{cor3}
For a semiring ${\mathcal S}=(S,+,\cdot,0,1)$ the following are equivalent:
\begin{enumerate}
\item[{\rm(i)}] ${\mathcal S}$ is a {\rm DNL}-semiring satisfying the identity $x\cdot n(n(y)+n(z))\approx n(n(x\cdot y)+n(x\cdot z))$.
\item[{\rm(ii)}] ${\bf L}({\mathcal S})=(S,\vee,\wedge,\cdot,\rightarrow,0,1)$ is a complete and prelinear residuated lattice satisfying the double negation law where $x\vee y=x+y$, $x\wedge y=n(n(x)+n(y))$, $x\odot y=x\cdot y$ and $x\rightarrow y=n(x\cdot n(y))$ for all $x,y\in S$.
\end{enumerate}
\end{Corollary}

\begin{proof}
According to Theorems~\ref{th4}, \ref{th6} and \ref{th7}, ${\mathcal S}$ is a DNL-semiring if and only if ${\bf L}({\mathcal S})$ is a complete residuated lattice satisfying the double negation law. Because of Theorem~2.40 in \cite B, ${\bf L}({\mathcal S})$ is prelinear if and only if it satisfies the identity
\[
x\cdot(y\wedge z)\approx(x\cdot y)\wedge(x\cdot z)
\]
which is equivalent to
\[
x\cdot n(n(y)+n(z))\approx n(n(x\cdot y)+n(x\cdot z)).
\]
\end{proof}

Authors' addresses:

Ivan Chajda \\
Palack\'y University Olomouc \\
Faculty of Science \\
Department of Algebra and Geometry \\
17.\ listopadu 12 \\
771 46 Olomouc \\
Czech Republic \\
ivan.chajda@upol.cz

Helmut L\"anger \\
TU Wien \\
Faculty of Mathematics and Geoinformation \\
Institute of Discrete Mathematics and Geometry \\
Wiedner Hauptstra\ss e 8-10 \\
1040 Vienna \\
Austria \\
helmut.laenger@tuwien.ac.at
\end{document}